\documentclass[a4paper]{amsart}
\usepackage{graphicx}
\usepackage{amsmath}
\usepackage{amssymb}
\usepackage{oldgerm}
\usepackage[all]{xy}

\newcommand{\Hom}{\operatorname{Hom}\nolimits}

\newcommand{\Pd}{\operatorname{pd}\nolimits}

\newcommand{\depth}{\operatorname{depth}\nolimits}

\newcommand{\Tor}{\operatorname{Tor}\nolimits}
\newcommand{\Ext}{\operatorname{Ext}\nolimits}

\newcommand{\m}{\operatorname{\mathfrak{m}}\nolimits}

\newcommand{\cx}{\operatorname{cx}\nolimits}

\newtheorem{theorem}{Theorem}[section]

\newtheorem*{condition}{Condition}

\newtheorem{lemma}[theorem]{Lemma}

\theoremstyle{definition}
\newtheorem*{example}{Example}

\theoremstyle{definition}

\theoremstyle{remark}

\theoremstyle{definition}

\begin{document}
\title{On the vanishing of (co)homology over local rings}
\author{Petter Andreas Bergh}
\address{Petter Andreas Bergh \newline Mathematical Institute \\ 24-29
  St.\ Giles \\ Oxford \\ OX1 3LB \\ United Kingdom} 
\email{bergh@maths.ox.ac.uk}
%\date{\today}

\thanks{The author acknowledges the support of the Marie Curie
Research Training Network through a LieGrits fellowship} 



\maketitle

\begin{abstract}
Considering modules of finite complete intersection dimension over
commutative Noetherian local rings, we prove (co)homology vanishing
results in which we assume the vanishing of \emph{nonconsecutive}
(co)homology groups. In fact, the (co)homology groups assumed to
vanish may be arbitrarily far apart from each other.
\end{abstract}

\section{Introduction}

The study of the rigidity of $\Tor$ for modules over commutative
Noetherian local rings began with Auslander's 1961 paper
\cite{Auslander}, which focused on torsion properties of tensor
products over unramified regular local rings. In this paper
Auslander proved his famous rigidity theorem; for a given pair of
finitely generated modules over an unramified regular local ring, if
one $\Tor$ group vanishes, then all subsequent $\Tor$ groups vanish.
This result was extended to all regular local rings by
Lichtenbaum in \cite{Lichtenbaum}, where the ramified case was
proved. Then in \cite{Murthy} Murthy generalized the rigidity
theorem to arbitrary complete intersections, by showing that if
$c+1$ consecutive $\Tor$ groups vanish for two finitely generated
modules over a complete intersection of codimension $c$, then all
the subsequent $\Tor$ groups also vanish. Three decades later Huneke
and Wiegand proved in \cite{Huneke} that given certain length and
dimension restrictions on the two modules involved, the vanishing
interval may be reduced by one, i.e.\ one only has to assume that
$c$ consecutive $\Tor$ groups vanish.

In \cite{Jorgensen1}, rather than considering vanishing intervals of
lengths determined by the codimension of the complete intersection,
Jorgensen used the notion of the \emph{complexity} of a module. He
proved that if $X$ is a finitely generated module of complexity $c$
over a complete intersection $R$, and
$$\Tor_n^R(X,Y) = \Tor_{n+1}^R(X,Y) = \cdots = \Tor_{n+c}^R(X,Y)=0$$
for a finitely generated module $Y$ and some integer $n > \dim R -
\depth X$, then $\Tor_i^R(X,Y)=0$ for all $i > \dim R - \depth X$.
In addition, by imposing the same length and dimension restrictions
on $X$ and $Y$ as Huneke and Wiegand did, he showed that the
vanishing interval also in this case can be reduced by one. These
results drastically generalize those of Murthy, Huneke and Wiegand,
since the complexity of a finitely generated module over a complete
intersection never exceeds the codimension of the ring. Moreover,
given a complete intersection $R$ of codimension $c>0$ and an
integer $t \in \{ 0,1, \dots, c \}$, there exist many finitely
generated $R$-modules having complexity $t$ (see, for example,
\cite[3.1-3.3]{Avramov3} and also \cite[Corollary 2.3]{Bergh2}). In
particular, there are a lot of modules of infinite projective
dimension having complexity strictly less than the codimension of
$R$. 

A fundamental feature shared by all these vanishing results is the
assumption that a certain number of \emph{consecutive} $\Tor$ groups
vanish, although some results assuming the vanishing of two families
of consecutive $\Tor$ groups were given in \cite{Jorgensen1} and
\cite{Jorgensen2}. The reason behind this is the change of rings
spectral sequence
$$\Tor_p^{R/(x)} \left ( X, \Tor_q^R(Y,R/(x)) \right )
\underset{p}{\Rightarrow} \Tor_{p+q}^R(X,Y)$$ associated to
$R/(x)$-modules $X$ and $Y$ (see \cite[XVI,{\S}5]{Cartan}), a
sequence which, when $x$ is a non-zerodivisor in $R$, degenerates
into a long exact sequence connecting $\Tor^{R/(x)}_i(X,Y)$ and
$\Tor^{R}_i(X,Y)$. In this situation the vanishing of $t$
consecutive $\Tor$ groups over $R/(x)$ implies the vanishing of
$t-1$ consecutive $\Tor$ groups over $R$, hence when the ring one
studies is a complete intersection one may use induction on the
codimension and the already known results for modules over
hypersurfaces and regular local rings. However, by using induction
on the complexity of a module rather than the codimension of the
ring, we give in this paper vanishing results for $\Tor$
(respectively, $\Ext$) in which we assume the vanishing of
\emph{nonconsecutive} $\Tor$ groups (respectively, $\Ext$ groups).
In fact, the (co)homology groups assumed to vanish may be
arbitrarily far apart from each other. This is easy to see in the
complexity one case; by a result of Eisenbud (see \cite{Eisenbud})
a finitely generated module of complexity one over a complete
intersection is eventually syzygy-periodic of period not more than
two, hence if one odd and one even (high enough) $\Tor$
(respectively, $\Ext$) vanish, then all higher $\Tor$ groups
(respectively, $\Ext$ groups) vanish.

As the proofs in this paper use induction on the complexity of a
module, we need to somehow be able to reduce to a situation in which
the module involved has complexity one less than the module we
started with. A method for doing this was studied in \cite{Bergh},
where it was shown that over a commutative Noetherian local ring
every module of finite complete intersection dimension has so-called
``reducible complexity". Therefore we state the main results for
modules having finite complete intersection dimension, a class of
modules containing all modules over complete intersections.

\section{Preliminaries}

Throughout we let ($A, \m, k$) be a local (meaning also commutative
Noetherian) ring, and we suppose all modules are finitely generated.
For an $A$-module $M$ with minimal free resolution
$$\bold{F}_M \colon \cdots \to F_2 \to F_1
\to F_0 \to M \to 0,$$ the rank of $F_n$, i.e.\ the integer $\dim_k
\Ext_A^n (M,k)$, is the $n$th \emph{Betti number} of $M$, and we
denote this by $\beta_n(M)$. The \emph{complexity} of $M$, denoted
$\cx M$, is defined as
$$\cx M = \inf \{ t \in \mathbb{N} \cup \{ 0
\} \mid \exists a \in \mathbb{R} \text{ such that } \beta_n(M) \leq
an^{t-1} \text{ for } n \gg 0 \}.$$ In general, the complexity of a
module may be infinite, whereas it is zero if and only if the module
has finite projective dimension. The $n$th \emph{syzygy} of $M$,
denoted $\Omega_A^n (M)$, is the cokernel of the map $F_{n+1} \to
F_n$, and it is unique up to isomorphism (note that $\Omega_A^0
(M)=M$ and $\cx M = \cx \Omega_A^i (M)$ for every $i \geq 0$). Now
let $N$ be an $A$-module, and consider a homogeneous element $\eta
\in \Ext^*_A(M,N)$. By choosing a map $f_{\eta} \colon
\Omega_A^{|\eta|}(M) \to N$ representing $\eta$, we obtain a
commutative pushout diagram
$$\xymatrix{
0 \ar[r] & \Omega_A^{|\eta|}(M) \ar[r] \ar[d]^{f_{\eta}} &
F_{|\eta|-1} \ar[r] \ar[d] & \Omega_A^{|\eta|-1}(M) \ar[r]
\ar@{=}[d] & 0 \\
0 \ar[r] & N \ar[r] & K_{\eta} \ar[r] & \Omega_A^{|\eta|-1}(M)
\ar[r] & 0 }$$ with exact rows. Note that the module $K_{\eta}$ is
independent, up to isomorphism, of the map $f_{\eta}$ chosen as a
representative for $\eta$. Following \cite{Bergh}, we say that $M$
has \emph{reducible complexity} if either the projective dimension
of $M$ is finite, or if the complexity of $M$ is positive and
finite, and there exists a homogeneous element $\eta \in
\Ext_A^*(M,M)$ such that $\cx K_{\eta} < \cx M$, $\depth K_{\eta} =
\depth M$, and $K_{\eta}$ has reducible complexity. The
cohomological element $\eta$ is said to \emph{reduce the complexity}
of $M$.

The module $M$ has finite \emph{complete intersection dimension} if
there exist local rings $R$ and $Q$ and a diagram $A \to R
\twoheadleftarrow Q$ of local homomorphisms such that $A \to R$ is
faithfully flat, $R \twoheadleftarrow Q$ is surjective with kernel
generated by a regular sequence (such a diagram $A \to R
\twoheadleftarrow Q$ is called a \emph{quasi-deformation} of $A$),
and $\Pd_Q (R \otimes_A M)$ is finite. Such modules were first
studied in \cite{Avramov3}, and the concept generalizes that of
virtual projective dimension defined in \cite{Avramov1}. As the name
suggests, modules having finite complete intersection dimension to a
large extent behave homologically like modules over complete
intersections. Indeed, over a complete intersection ($S,
\mathfrak{n}$) \emph{every} module has finite complete intersection
dimension; the completion $\widehat{S}$ of $S$ with respect to the
$\mathfrak{n}$-adic topology is the residue ring of a regular local
ring $Q$ modulo an ideal generated by a regular sequence, and so $S
\to \widehat{S} \twoheadleftarrow Q$ is a quasi deformation.

By \cite[Proposition 2.2]{Bergh} every module of finite complete
intersection dimension has reducible complexity, and given such a
module the reducing process decreases the complexity by exactly one.
Moreover, the following result shows that by passing to a suitable
faithfully flat extension of the ring, we may assume that the
cohomological element reducing the complexity of the module is of
degree two. Moreover, if $\eta$ is a cohomological element reducing
the complexity of the module, then every positive power of $\eta$
also reduces the complexity.

\begin{lemma}\label{reducing}
Let $M$ be an $A$-module of finite complete intersection dimension.
\begin{enumerate}
\item[(i)] There exists a quasi-deformation $A \to R \twoheadleftarrow Q$
such that the $R$-module $R \otimes_A M$ has reducible complexity by
an element $\eta \in \Ext_R^2(R \otimes_A M,R \otimes_A M)$.
\item[(ii)] If $\eta \in \Ext_A^*(M,M)$ reduces the complexity of
$M$, then so does $\eta^t$ for every $t \geq 1$.
\end{enumerate}
\end{lemma}

\begin{proof}
(i) By \cite[Proposition 7.2(2)]{Avramov3} there exists a
quasi-deformation $A \to R \twoheadleftarrow Q$ in which the
kernel of the map $R \twoheadleftarrow Q$ is generated by a
$Q$-regular element, and such that the Eisenbud operator (with
respect to $Q$) on the minimal free resolution of the $R$-module
$R \otimes_A M$ is eventually surjective. This operator is a chain
endomorphism corresponding to an element $\eta \in \Ext_R^2(R
\otimes_A M,R \otimes_A M)$, and the proof of \cite[Proposition
2.2(i)]{Bergh} shows that $\cx_R K_{\eta} = \cx_R (R \otimes_A
M)-1$ and $\depth_R K_{\eta} = \depth_R (R \otimes_A M)$. It also
follows from the proof of \cite[Proposition 7.2(2)]{Avramov3} that
there exists a local ring $\overline{Q}$ and an epimorphism $R
\twoheadleftarrow \overline{Q}$ factoring through $R
\twoheadleftarrow Q$, whose kernel is generated by a
$\overline{Q}$-regular sequence, and such that
$\Pd_{\overline{Q}}(R \otimes_A M)$ is finite. Since
$\Pd_{\overline{Q}}R < \infty$ and $R$ is a flat $A$-module, we
see that $\Pd_{\overline{Q}} \left ( \Omega_R^1(R \otimes_A M)
\right )$ is finite, and the exact sequence
$$0 \to R \otimes_A M \to K_{\eta} \to \Omega_R^1(R \otimes_A
M) \to 0$$ shows that the same is true for $\Pd_{\overline{Q}}
K_{\eta}$. Therefore $K_{\eta}$ has finite complete intersection
dimension and consequently has reducible complexity, and this shows
that the $R$-module $R \otimes_A M$ has reducible complexity.

(ii) Suppose $\eta^t$ reduces the complexity of $M$ for some $t
\geq 1$, and let $A \to R \twoheadleftarrow Q$ be a quasi
deformation for which $\Pd_Q(R \otimes_A M)$ is finite. As $\Pd_Q
R$ is finite, so is $\Pd_Q \left ( R \otimes_A \Omega_A^i(M)
\right )$ for any $i \geq 0$, hence from the exact sequence
$$0 \to R \otimes_A M \to R \otimes_A K_{\eta^t} \to R \otimes_A
\Omega_A^{t|\eta|-1}(M) \to 0$$ we see that $\Pd_Q (R \otimes_A
K_{\eta^t}) < \infty$. Now from \cite[Lemma 2.3]{Bergh} there
exists an exact sequence
$$0 \to \Omega_A^{|\eta|}(K_{\eta^t}) \to K_{\eta^{t+1}} \oplus F
\to K_{\eta} \to 0$$ in which $F$ is a free $A$-module, and
tensoring this sequence with $R$ we see that $\Pd_Q (R \otimes_A
K_{\eta^{t+1}} )$ is finite. Consequently the $A$-module
$K_{\eta^{t+1}}$ has finite complete intersection dimension, in
particular its complexity is reducible. Moreover, the above exact
sequence gives
$$\cx K_{\eta^{t+1}} \leq \max \{ \cx K_{\eta}, \cx K_{\eta^t} \}
< \cx M,$$ and this shows that $\eta^{t+1}$ reduces the complexity
of $M$.
\end{proof}

We end this section with recalling the following facts regarding
flat extensions of local rings. Let $S \to T$ be a faithfully flat
local homomorphism, and let $X$ and $Y$ be $S$-modules. If
$\bold{F}_X$ is a minimal $S$-free resolution of $X$, then the
complex $T \otimes_S \bold{F}_X$ is a minimal $T$-free resolution of
$T \otimes_S X$, and by \cite[Proposition (2.5.8)]{Grothendieck}
there exist natural isomorphisms
\begin{eqnarray*}
\Hom_T (T \otimes_S \bold{F}_X, T \otimes_S Y) & \simeq & T
\otimes_S \Hom_S ( \bold{F}_X, Y), \\
(T \otimes_S \bold{F}_X) \otimes_T (T \otimes_S Y) & \simeq & T
\otimes_S ( \bold{F}_X \otimes_S Y).
\end{eqnarray*}
This gives isomorphisms
\begin{eqnarray*}
\Ext_T^i (T \otimes_S X,T \otimes_S Y) & \simeq & T \otimes_S
\Ext_S^i (X,Y), \\
\Tor^T_i (T \otimes_S X,T \otimes_S Y) & \simeq & T \otimes_S
\Tor^S_i (X,Y),
\end{eqnarray*}
and as $T$ is faithfully $S$-flat we then get
\begin{eqnarray*}
\Ext_T^i (T \otimes_S X,T \otimes_S Y)=0 & \Leftrightarrow &
\Ext_S^i (X,Y) =0, \\
\Tor^T_i (T \otimes_S X,T \otimes_S Y)=0 & \Leftrightarrow &
\Tor^S_i (X,Y) =0.
\end{eqnarray*}
Moreover, there are equalities
\begin{eqnarray*}
\cx_S X &=& \cx_T (T \otimes_S X) \\
\depth_S X - \depth_S Y &=& \depth _T (T \otimes_S X) - \depth_T (T
\otimes_S Y),
\end{eqnarray*}
where the one involving depth follows from \cite[Theorem
23.3]{Matsumura}.

\section{Vanishing results}

Throughout this section we fix two nonzero $A$-modules $M$ and
$N$, with $M$ of complexity $c$. The first main result shows that
the vanishing of $\Ext$ for a certain sequence of numbers forces
the vanishing of all the higher $\Ext$ groups. The integers for
which the cohomology groups are assumed to vanish may be
arbitrarily far apart from each other.

\begin{theorem}\label{cohomologyvanishing}
Suppose $M$ has finite complete intersection dimension. If there
exist an odd number $q \geq 1$ and a number $n> \depth A - \depth
M$ such that $\Ext^i_A(M,N)=0$ for $i \in \{ n, n+q, \dots, n+cq
\}$, then $\Ext^i_A(M,N)=0$ for all $i> \depth A - \depth M$.
\end{theorem}

\begin{proof}
We prove this by induction on $c$. If $c =0$, then the projective
dimension of $M$ is finite and equal to $\depth A - \depth M$ by
the Auslander-Buchsbaum formula, and so the theorem holds in this
case. Now suppose $c$ is positive, and write $q$ as $q=2t-1$ where
$t \geq 1$. By Lemma \ref{reducing}(i) there exists a quasi
deformation $A \to R \twoheadleftarrow Q$ such that the $R$-module
$R \otimes_A M$ has reducible complexity by an element $\eta \in
\Ext_R^2(R \otimes_A M,R \otimes_A M)$. Moreover, as $R \otimes_A
M$ has finite complete intersection dimension we see from Lemma
\ref{reducing}(ii) that the element $\eta^t$ also reduces the
complexity of $R \otimes_A M$.

The exact sequence
\begin{equation*}\label{ES}
0 \to R \otimes_A M \to K_{\eta^t} \to \Omega_R^q(R \otimes_A M) \to
0 \tag{$\dagger$}
\end{equation*}
induces a long exact sequence
$$\cdots \to \Ext_R^{i+q}({_RM},{_RN}) \to \Ext_R^i(K_{\eta^t},{_RN})
\to \Ext_R^i({_RM},{_RN}) \to \cdots$$ of cohomology groups (where,
for an $A$-module $X$, the notation ${_RX}$ is short hand for $R
\otimes_A X$), and so the vanishing of the $c+1$ cohomology groups
$\Ext_A(M,N)$ implies that $\Ext_R^i(K_{\eta^t}, R \otimes_A N) =0$
for $i \in \{ n, n+q, \dots, n+(c-1)q \}$. As in the proof of Lemma
\ref{reducing} the $R$-module $K_{\eta^t}$ has finite complete
intersection dimension, and the equality $\depth_R K_{\eta^t} =
\depth_R (R \otimes_A M)$ holds. Therefore, since $\depth_A A -
\depth_A M = \depth_R R - \depth_R (R \otimes_A M)$ and $\cx
K_{\eta^t} = \cx (R \otimes_A M)-1 = c-1$, induction gives
$\Ext_R^i(K_{\eta^t}, R \otimes_A N) =0$ for $i > \depth_A A -
\depth_A M$. Consequently
$$\Ext_R^i(R \otimes_A M, R \otimes_A N) \simeq \Ext_R^{i+j(q+1)}(R
\otimes_A M, R \otimes_A N)$$ for $i > \depth_A A - \depth_A M$ and
$j \geq 0$, and by considering all the pairs $(i,j) \in \{ (n,c),
(n+q, c-1), \dots, (n+cq,0) \}$ we see that $\Ext^i_R(R \otimes_A M,
R \otimes_A N)$ vanishes for $n+cq \leq i \leq n+cq+c$. Then
$\Ext^i_A(M,N)$ also vanishes for these $c+1$ consecutive values,
and by \cite[Corollary 2.3]{Jorgensen2} and \cite[Theorem
1.4]{Avramov3} we get $\Ext^i_A(M,N)=0$ for $i > \depth_A A -
\depth_A M$.
\end{proof}

The next result is a homology version of Theorem
\ref{cohomologyvanishing}, and the proof is similar.

\begin{theorem}\label{homologyvanishing}
Suppose $M$ has finite complete intersection dimension. If there
exist an odd number $q \geq 1$ and a number $n> \depth A - \depth
M$ such that $\Tor^A_i(M,N)=0$ for $i \in \{ n, n+q, \dots, n+cq
\}$, then $\Tor^A_i(M,N)=0$ for all $i> \depth A - \depth M$.
\end{theorem}

\begin{proof}
By induction on $c$, where the case $c=0$ follows from the
Auslander-Buchsbaum formula. Suppose therefore $c$ is positive,
and let the notation be as in the proof of Theorem
\ref{cohomologyvanishing}. The exact sequence (\ref{ES}) induces a
long exact sequence
$$\cdots \to \Tor_i^R ({_RM},{_RN}) \to \Tor_i^R
(K_{\eta^t},{_RN}) \to \Tor_{i+q}^R ({_RM},{_RN}) \to \cdots$$ of
homology modules, giving $\Tor^R_i(K_{\eta^t},{_RN})=0$ for $i \in
\{ n, n+q, \dots, n+(c-1)q \}$. By induction
$\Tor^R_i(K_{\eta^t},{_RN})$ vanishes for all $i > \depth_A A -
\depth_A M$, giving an isomorphism
$$\Tor^R_i(R \otimes_A M, R \otimes_A N) \simeq \Tor^R_{i+j(q+1)}(R
\otimes_A M, R \otimes_A N)$$ for $i > \depth_A A - \depth_A M$ and
$j \geq 0$. As in the previous proof we see that $\Tor^R_i(R
\otimes_A M, R \otimes_A N)$ vanishes for $n+cq \leq i \leq n+cq+c$,
and therefore $\Tor^A_i(M,N)$ also vanishes for these $c+1$
consecutive values. From \cite[Corollary 2.6]{Jorgensen2} and
\cite[Theorem 1.4]{Avramov3} we see that $\Tor_i^A(M,N)=0$ for $i >
\depth_A A - \depth_A M$.
\end{proof}

The number of homology and cohomology groups assumed to vanish in
the above two theorems cannot be reduced in general; there exist
examples of a complete intersection $R$ and modules $X$ and
$Y$ for which $\cx X$ consecutive values of $\Tor^R(X,Y)$ or
$\Ext_R(X,Y)$ vanish (starting beyond $\dim R$), while not all the
higher homology or cohomology groups vanish (see
\cite[4.1]{Jorgensen1} for a homology example). We include such an
example, which also shows that the distance between the
(co)homology groups assumed to vanish cannot be an even integer.

\begin{example}
Let $k$ be a field, let $k[[X,Y]]$ be the ring of formal power
series in the indeterminates $X$ and $Y$, and denote the one
dimensional hypersurface $k[[X,Y]]/(XY)$ by $A$. The complex
$$\cdots \to A \xrightarrow{x} A \xrightarrow{y} A \xrightarrow{x} A
\to A/(x) \to 0$$ is exact, and is therefore a minimal free
resolution of the module $M = A/(x)$. Denoting by $N$ the module
$A/(y)$, we see that $\Tor_i^A(M,N)$ vanishes when $i$ is odd and is
nonzero when $i$ is non-negative and even, and that $\Ext^i_A(M,N)$
vanishes when $i$ is even and is nonzero when $i$ is positive and
odd.
\end{example}

However, when the ring is a complete intersection and one of the
modules for which we are computing (co)homology has finite length,
we may reduce the number of (co)homology groups assumed to vanish by
one. To prove this, we need the following lemma, which is really
just a special case of the two theorems to be proven next. Namely,
it treats the case when the (co)homology groups assumed to vanish
are all consecutive. Note that the homology part of the lemma
follows easily from \cite[Theorem 2.6]{Jorgensen1}, whereas the
complexity one case of the cohomology part, for which we include a
proof, is more or less similar to the proof of \cite[Proposition
2.2]{Huneke}.

\begin{lemma}\label{finitelength}
Let $A$ be a complete intersection, and suppose $N$ has finite
length.
\begin{enumerate}
\item[(i)] If there exists a number $n> \dim A - \depth M$ such
that $\Ext^i_A(M,N)=0$ for $n \leq i \leq n+c-1$, then
$\Ext^i_A(M,N)=0$ for all $i> \dim A - \depth M$. \item[(ii)] If
there exists a number $n> \dim A - \depth M$ such that
$\Tor^A_i(M,N)=0$ for $n \leq i \leq n+c-1$, then
$\Tor^A_i(M,N)=0$ for all $i> \dim A - \depth M$.
\end{enumerate}
\end{lemma}

\begin{proof}
(i) As length is preserved under faithfully flat extensions, we
may assume the ring $A$ is complete and has infinite residue field
(if the latter is not the case, then we use the faithfully flat
extension $A \to A[x]_{\m A[x]}$, where $x$ is an indeterminate).
We argue by induction on $c$, where the case $c=0$ follows from
the Auslander-Buchsbaum formula. Suppose therefore $c=1$, and
choose, by \cite[Theorem 1.3]{Jorgensen1}, a local ring $R$ and a
non-zerodivisor $x \in R$ such that $A =R/(x)$ and $\Pd_R M <
\infty$. As $\dim R = \dim A+1$ and $\depth_R M = \depth_A M$, the
Auslander-Buchsbaum formula gives $\Pd_R M = \dim A - \depth_A M
+1$.

\sloppy The change of rings spectral sequence (see
\cite[XVI,{\S}5]{Cartan})
$$\Ext_A^p \left ( M, \Ext^q_R(A,N) \right )
\underset{p}{\Rightarrow} \Ext^{p+q}_R(M,N)$$ degenerates into a
long exact sequence
$$\xymatrix@R=0.2pc@C=1pc{
0 \ar[r] & \Ext_A^1(M,N) \ar[r] & \Ext_R^1(M,N) \ar[r] & \Hom_A(M,N)
\\
\ar[r] & \Ext_A^2(M,N) \ar[r] & \Ext_R^2(M,N) \ar[r] & \Ext_A^1(M,N)
\\
& \vdots & \vdots & \vdots \\
\ar[r] & \Ext_A^i(M,N) \ar[r] & \Ext_R^i(M,N) \ar[r] &
\Ext_A^{i-1}(M,N) \\
& \vdots & \vdots & \vdots }$$ connecting the cohomology groups
over $A$ to those over $R$. Denote $\Pd_R M$ by $d$, and by
$\beta_m$ the $m$th Betti number of $M$ over $R$. Localizing the
minimal $R$-free resolution
$$0 \to R^{\beta_d} \to \cdots \to R^{\beta_1} \to R^{\beta_0} \to
M \to 0$$ with respect to the multiplicatively closed set $\{ 1,
x, x^2, \dots \} \subseteq R$, keeping in mind that $xM=0$, we
obtain an exact sequence
$$0 \to R_x^{\beta_d} \to \cdots \to R_x^{\beta_1} \to
R_x^{\beta_0} \to 0$$ from which we conclude that the equality
$\sum_{i=0}^d(-1)^i\beta_i =0$ holds. As $\beta_m = \ell \left (
\Ext_R^m(M,k) \right )$, where $\ell (-)$ denotes the length
function, we see that $\sum_{i=0}^d(-1)^i \ell \left (
\Ext_R^i(M,k) \right ) =0$, and then by induction on length
$\sum_{i=0}^d(-1)^i \ell \left ( \Ext_R^i(M,X) \right ) =0$ for
every $R$-module $X$ of finite length. Now let $m > \dim A -
\depth_A M$ be an integer. By taking the alternate sum of the
lengths in the above long exact sequence, ending in
$\Ext_R^{m+1}(M,N)=0$ and using the isomorphism $\Hom_A(M,N)
\simeq \Hom_R(M,N)$, we see that the only terms not canceling are
$\ell \left ( \Ext_A^m(M,N) \right )$ and $\ell \left (
\Ext_A^{m+1}(M,N) \right )$. As this alternate sum is zero we get
$\ell \left ( \Ext_A^m(M,N) \right ) = \ell \left (
\Ext_A^{m+1}(M,N) \right )$, and it follows immediately that if
$\Ext_A^n(M,N)=0$ for some $n > \dim A - \depth_A M$, then
$\Ext_A^i(M,N)$ vanishes for all $i > \dim A - \depth_A M$. This
proves the case $c=1$.

Now suppose $c>1$. As $A$ is complete there exists an element
$\eta \in \Ext_A^2(M,M)$ reducing the complexity of $M$. The exact
sequence
$$0 \to M \to K_{\eta} \to \Omega_A^1(M) \to 0$$
induces a long exact sequence
$$\cdots \to \Ext_A^{i+1}(M,N) \to \Ext_A^i(K_{\eta},N) \to
\Ext_A^i(M,N) \to \cdots$$ in $\Ext$, giving
$\Ext_A^i(K_{\eta},N)=0$ for $n \leq i \leq n+c-2$. By induction
$\Ext_A^i(K_{\eta},N)$ vanishes for all $n > \dim A - \depth M$,
and the long exact sequence then implies $\Ext_A^i(M,N)=0$ for all
$n > \dim A - \depth M$.

(ii) The homology part is a special case of \cite[Theorem
2.6]{Jorgensen1}. It can also be proved by an argument similar to
the above.
\end{proof}

We are now ready to prove the results which sharpen Theorem
\ref{cohomologyvanishing} and Theorem \ref{homologyvanishing} when
$N$ has finite length. Namely, it is enough to assume that $c$
(co)homology groups vanish.

\begin{theorem}\label{cohomologyvanishinglength}
Let $A$ be a complete intersection, and suppose $N$ has finite
length. If there exist an odd number $q \geq 1$ and a number $n>
\dim A - \depth M$ such that $\Ext^i_A(M,N)=0$ for $i \in \{ n, n+q,
\dots, n+(c-1)q \}$, then $\Ext^i_A(M,N)=0$ for all $i> \dim A -
\depth M$.
\end{theorem}

\begin{proof}
If $c=0$ then the vanishing assumption is vacuous, but as $M$ has
finite projective dimension the conclusion follows from the
Auslander-Buchsbaum formula. The case $c=1$ follows from Lemma
\ref{finitelength}(i), so suppose $c>1$ and that the result holds if
we replace $M$ by a module having complexity $c-1$. As before we may
suppose $A$ is complete, hence there exists an element $\eta \in
\Ext_A^2(M,M)$ reducing the complexity of $M$. Write $q$ as
$q=2t-1$, where $t \geq 1$. The element $\eta^t$ also reduces the
complexity of $M$, and the corresponding exact sequence
$$0 \to M \to K_{\eta^t} \to \Omega_A^q(M) \to 0$$
induces a long exact sequence
$$\cdots \to \Ext_A^{i+q}(M,N) \to \Ext_A^i(K_{\eta^t},N) \to
\Ext_A^i(M,N) \to \cdots$$ in $\Ext$. The assumption on the
vanishing $\Ext_A(M,N)$ groups gives $\Ext_A^i(K_{\eta^t},N)=0$ for
$i \in \{ n, n+q, \dots, n+(c-2)q \}$, hence by induction
$\Ext_A^i(K_{\eta^t},N)$ vanishes for all $n > \dim A - \depth M$.
But then for every $i > \dim A - \depth M$ and $j \geq 0$ there is
an isomorphism
$$\Ext_A^i(M,N) \simeq \Ext_A^{i+j(q+1)}(M,N),$$
and by considering the pairs $(i,j) \in \{ (n,c), (n+q,c-1), \dots,
(n+(c-1)q,1) \}$ we see that $\Ext_A^i(M,N)$ vanishes for $n+cq+1
\leq i \leq n+cq+c$. By Lemma \ref{finitelength}(i) we are done.
\end{proof}

We omit the proof of the homology version of Theorem
\ref{cohomologyvanishinglength}, as it is completely analogous to
the proofs of Theorem \ref{homologyvanishing} and Theorem
\ref{cohomologyvanishinglength}.

\begin{theorem}\label{homologyvanishinglength}
Let $A$ be a complete intersection, and suppose $N$ has finite
length. If there exist an odd number $q \geq 1$ and a number $n>
\dim A - \depth M$ such that $\Tor^A_i(M,N)=0$ for $i \in \{ n, n+q,
\dots, n+(c-1)q \}$, then $\Tor^A_i(M,N)=0$ for all $i> \dim A -
\depth M$.
\end{theorem}

An obvious question that arises from Theorem
\ref{cohomologyvanishing} and Theorem \ref{homologyvanishing} is
whether or not the gaps between the (co)homology groups assumed to
vanish have to be of the \emph{same} length. Consider therefore the
following cohomology vanishing condition:

\begin{condition}
There exist positive odd integers $q_1, \dots, q_c$ such that
$$\Ext_A^n(M,N) = \Ext_A^{n+q_1}(M,N) = \cdots = \Ext_A^{n+q_1+ \cdots
+q_c}(M,N)=0$$ for some $n> \depth A - \depth M$.
\end{condition}

Does the conclusion of Theorem \ref{cohomologyvanishing} hold if we
replace the vanishing assumption in the theorem with the
considerably weaker condition stated above? Similarly we may ask if
the conclusion of Theorem \ref{homologyvanishing} holds if we
replace the vanishing assumption in the theorem with the homology
version of the above condition.

In case the answer to the above questions are positive, it is
possible that a proof similar to those used in the main results
exist, i.e.\ a proof based on reducing the complexity of $M$. We end
this paper with two results showing that the answer is positive when
the complexity of $M$ is not more than $2$. A proof is provided only
for the cohomology version, the proof of the homology version is
analogous.

\begin{theorem}\label{cohomologytwo}
Suppose $M$ has finite complete intersection dimension and
complexity $2$. If there exist odd numbers $p \geq 1$ and $q \geq 1$
such that
$$\Ext_A^n(M,N) = \Ext_A^{n+p}(M,N)= \Ext_A^{n+p+q}(M,N)=0$$ for some
$n> \depth A - \depth M$, then $\Ext^i_A(M,N)=0$ for all $i> \depth
A - \depth M$.
\end{theorem}

\begin{proof}
We argue by induction on the sum $p+q$ (by assumption this sum is at
least $2$). If $p+q=2$, then $p=q=1$, and we are done by Theorem
\ref{cohomologyvanishing}. Suppose therefore $p+q>2$, and write $p$
and $q$ as $p=2s-1$ and $q=2t-1$, where $s,t \in \mathbb{N}$.

Choose, by Lemma \ref{reducing}(i), a quasi deformation $A \to R
\twoheadleftarrow Q$ such that the complexity of the $R$-module $R
\otimes_A M$ is reducible by an element $\eta \in \Ext_R^2(R
\otimes_A M,R \otimes_A M)$. For an $A$-module $X$, denote by $_RX$
the $R$-module $R \otimes_A X$. By Lemma \ref{reducing}(ii) the
element $\eta^s$ also reduces the complexity of ${_RM}$, and from
its associated exact sequence
$$0 \to {_RM} \to K_{\eta^s} \to \Omega_R^p({_RM})
\to 0$$ we deduce that $\Ext_R^n(K_{\eta^s}, {_RN})=0$. However, the
$R$-module $K_{\eta^s}$ has complexity $1$, and therefore, by
\cite[Theorem 7.3]{Avramov3}, the modules $\Omega_R^i(K_{\eta^s})$
and $\Omega_R^{i+2}(K_{\eta^s})$ are isomorphic for $i > \depth_R R
- \depth_R K_{\eta^s}$. Consequently, as $\depth_R K_{\eta^s} =
\depth_R ({_RM})$, we see that
$$\Ext_R^{n+2i}(K_{\eta^s}, {_RN})=0$$
for every $i \geq 0$, in particular $\Ext_R^{n+q-1}(K_{\eta^s},
{_RN})=0$. Now since the above exact sequence yields an exact
sequence
$$\Ext_R^{n+q-1}(K_{\eta^s}, {_RN}) \to \Ext_R^{n+q-1}({_RM},{_RN}) \to
\Ext_R^{n+q}(\Omega_R^p({_RM}),{_RN}),$$ whose end terms both
vanish, we see that $\Ext_R^{n+q-1}({_RM},{_RN})=0$.

The element $\eta^t \in \Ext_R^*({_RM},{_RM})$ also reduces the
complexity of ${_RM}$, and from its associated exact sequence
$$0 \to {_RM} \to K_{\eta^t} \to \Omega_R^q({_RM})
\to 0$$ we see that $\Ext_R^{n+p}(K_{\eta^t}, {_RN})=0$. As with the
module $K_{\eta^s}$, the module $K_{\eta^t}$ is eventually periodic
of period $2$, and therefore
$$\Ext_R^{n+1+2i}(K_{\eta^t}, {_RN})=0$$
for every $i \geq 0$. In particular $\Ext_R^{n+1}(K_{\eta^t},
{_RN})=0$, and since the previous exact sequence yields an exact
sequence
$$\Ext_R^n({_RM},{_RN}) \to \Ext_R^{n+1}(\Omega_R^q({_RM}),{_RN})
\to \Ext_R^{n+1}(K_{\eta^t}, {_RN}),$$ we see that
$\Ext_R^{n+q+1}({_RM},{_RN})=0$.

We have just shown that both $\Ext_R^{n+q-1}({_RM},{_RN})$ and
$\Ext_R^{n+q+1}({_RM},{_RN})$ vanish. Now since $p$ and $q$ are not
equal, either $p > q \geq 1$ or $q > p \geq 1$. Define integers
$n',p'$ and $q'$ by 
$$(n',p',q') = \left \{ \begin{array}{ll}
                 (n+q+1,p-q-1,q) & \text{if } p>q \\
                 (n,p,q-p-1) & \text{if } p<q,
                \end{array} \right. $$
Then $p'$ and $q'$ are both odd, the groups $\Ext_R^i({_RM},{_RN})$
vanish for $i \in \{ n', n'+p', n'+p'+q' \}$, and $p'+q' < p+q$. By
induction $\Ext_R^i({_RM},{_RN})$ vanishes for $i> \depth A - \depth
M$, and consequently $\Ext_A^i(M,N)=0$ for $i> \depth A - \depth M$.
\end{proof}

\begin{theorem}\label{homologytwo}
Suppose $M$ has finite complete intersection dimension and
complexity $2$. If there exist odd numbers $p \geq 1$ and $q \geq 1$
such that
$$\Tor^A_n(M,N) = \Tor^A_{n+p}(M,N)= \Tor^A_{n+p+q}(M,N)=0$$ for some
$n> \depth A - \depth M$, then $\Tor_i^A(M,N)=0$ for all $i> \depth
A - \depth M$.
\end{theorem}

\end{document}